\documentclass[11pt]{article} 
\usepackage{latexsym} 
\usepackage{amsmath}
\usepackage{amssymb} 
\usepackage{amscd} 
\usepackage{graphicx}
\begin{document} 
\newtheorem{Th}{Theorem}[section]
\newtheorem{Cor}{Corollary}[section]
\newtheorem{Prop}{Proposition}[section]
\newtheorem{Lem}{Lemma}[section]
\newtheorem{Def}{Definition}[section]
\newtheorem{Rem}{Remark}[section]
\newtheorem{Ex}{Example}[section]
\newtheorem{stw}{Proposition}[section]


\newcommand{\bet}{\begin{Th}}
\newcommand{\ent}{\stepcounter{Cor}
   \stepcounter{Prop}\stepcounter{Lem}\stepcounter{Def}
   \stepcounter{Rem}\stepcounter{Ex}\end{Th}}


\newcommand{\bec}{\begin{Cor}}
\newcommand{\enc}{\stepcounter{Th}
   \stepcounter{Prop}\stepcounter{Lem}\stepcounter{Def}
   \stepcounter{Rem}\stepcounter{Ex}\end{Cor}}
\newcommand{\bep}{\begin{Prop}}
\newcommand{\enp}{\stepcounter{Th}
   \stepcounter{Cor}\stepcounter{Lem}\stepcounter{Def}
   \stepcounter{Rem}\stepcounter{Ex}\end{Prop}}
\newcommand{\bel}{\begin{Lem}}
\newcommand{\enl}{\stepcounter{Th}
   \stepcounter{Cor}\stepcounter{Prop}\stepcounter{Def}
   \stepcounter{Rem}\stepcounter{Ex}\end{Lem}}
\newcommand{\bef}{\begin{Def}}
\newcommand{\enf}{\stepcounter{Th}
   \stepcounter{Cor}\stepcounter{Prop}\stepcounter{Lem}
   \stepcounter{Rem}\stepcounter{Ex}\end{Def}}
\newcommand{\ber}{\begin{Rem}}
\newcommand{\enr}{
   \stepcounter{Th}\stepcounter{Cor}\stepcounter{Prop}
   \stepcounter{Lem}\stepcounter{Def}\stepcounter{Ex}\end{Rem}}
\newcommand{\bee}{\begin{Ex}}
\newcommand{\ene}{
   \stepcounter{Th}\stepcounter{Cor}\stepcounter{Prop}
   \stepcounter{Lem}\stepcounter{Def}\stepcounter{Rem}\end{Ex}}
\newcommand{\Proof}{\noindent{\it Proof\,}:\ }

\newcommand{\EE}{\mathbf{E}}
\newcommand{\QQ}{\mathbf{Q}}
\newcommand{\R}{\mathbf{R}}
\newcommand{\C}{\mathbf{C}}
\newcommand{\ZZ}{\mathbf{Z}}
\newcommand{\NN}{\mathbf{N}}
\newcommand{\PP}{\mathbf{P}}
\newcommand{\uuu}{{u}}
\newcommand{\xxx}{{x}}
\newcommand{\aaa}{{a}}
\newcommand{\bbb}{{b}}
\newcommand{\AAA}{\mathbf{A}}
\newcommand{\BBB}{\mathbf{B}}
\newcommand{\ccc}{{c}}
\newcommand{\iii}{{i}}
\newcommand{\jjj}{{j}}
\newcommand{\kkk}{{k}}
\newcommand{\rrr}{{r}}
\newcommand{\FFF}{{F}}
\newcommand{\yyy}{{y}}
\newcommand{\ppp}{{p}}
\newcommand{\qqq}{{q}}
\newcommand{\nnn}{{n}}
\newcommand{\vvv}{{v}}
\newcommand{\eee}{{e}}
\newcommand{\fff}{{f}}
\newcommand{\www}{{w}}
\newcommand{\0}{{0}}
\newcommand{\lon}{\longrightarrow}
\newcommand{\ga}{\gamma}
\newcommand{\pa}{\partial}
\newcommand{\QED}{\hfill $\Box$}
\newcommand{\id}{{\mbox {\rm id}}}
\newcommand{\Ker}{{\mbox {\rm Ker}}}
\newcommand{\Tan}{{\mbox {\rm Tan}}}
\newcommand{\grad}{{\mbox {\rm grad}}}
\newcommand{\ind}{{\mbox {\rm ind}}}
\newcommand{\rot}{{\mbox {\rm rot}}}
\newcommand{\diver}{{\mbox {\rm div}}}
\newcommand{\Gr}{{\mbox {\rm Gr}}}
\newcommand{\rank}{{\mbox {\rm rank}}}
\newcommand{\ord}{{\mbox {\rm ord}}}
\newcommand{\sign}{{\mbox {\rm sign}}}
\newcommand{\Spin}{{\mbox {\rm Spin}}}
\newcommand{\Symp}{{\mbox {\rm Sp}}}
\newcommand{\Int}{{\mbox {\rm Int}}}
\newcommand{\codim}{{\mbox {\rm codim}}}
\def\mod{{\mbox {\rm mod}}}
\newcommand{\qed}{\hfill $\Box$ \par}


\title{Normal and tangent maps to frontals} 

\author{Goo \textsc{Ishikawa}\thanks{Department of Mathematics, Faculty of Sciences, Hokkaido University, Sapporo 060-0810, Japan. 
e-mail: 
ishikawa@math.sci.hokudai.ac.jp}
}

\renewcommand{\thefootnote}{\fnsymbol{footnote}}
\footnotetext{Supported by JSPS and RFBR under the Japan-Russia Research Cooperative Program 120194801 and 
JSPS KAKENHI Grant Number JP19K03458.}

\date{ }

\maketitle

\begin{abstract}
\noindent
The notion of frontals in Euclidean space is introduced and the normal and tangent maps to 
frontals are studied for both geometrical and dynamical aspects of frontals. 
Moreover we observe that parallels of the tangent map to a frontal curve is right equivalent to the tangent map of a frontal curve under some natural conditions. 
\end{abstract}

\section{Introduction.} 

Given a submanifold in Euclidean space, we consider 
two typical mappings: 
the normal map which is ruled by normal spaces along the submanifold,  
and the tangent map, ruled by tangent spaces. 

Let $f : U^n \hookrightarrow \R^{n+p}$ be an $n$-dimensional submanifold in the Euclidean 
space $\R^{n+p}$. 
Then the {\it normal map} 
${\mathrm{Nor}}(f) : NU \to \R^{n+p}$ from the normal bundle $NU$ of $U$ is defined by 
$
{\mathrm{Nor}}(f)(t, \nu) := t + \nu, \ 
$ 
for $t \in U$ and $\nu \in N_tU$, using the affine structure of $\R^{n+p}$. 
Normal maps naturally appear, for example, 
in geometric optics and extrinsic differential geometry (\cite{Arnold, PT}). 
Their singularities are called {\it caustics} or {\it focal set}. 
Singularities of normal maps are regarded as Lagrangian singularities. 
They are studied and classified in Lagrangian singularity theory, 
by, for instance, the method of generating families (\cite{AGV, AGLV, Bruce}). 
Moreover singularities of canal hypersurfaces of submanifolds and 
parallels of hypersurfaces (resp. curves), for instance, are studied via 
the restrictions of normal maps by Legendrian singularity theory as 
singularities of {\it wavefronts} (\cite{AGV, Arnold}). Note that 
Legendrian singularity theory has wide applications also to the study of 
differential equations (see \cite{DIIS} for instance). 

The {\it tangent map} 
$
\Tan(f) : TU \to \R^{n+p}
$ 
of $f$ from the tangent bundle $TU$ is defined by 
$
\Tan(f)(t, \tau) := t + \tau, \ 
$ 
for $t \in U$ and $\tau \in T_tU$. 
Contrary to normal maps, tangent maps of submanifolds, 
which are natural subjects as well, 
have not necessarily Lagrangian singularities,  
and tangent maps have very 
degenerate singularities even in generic cases. 
For example, tangent developables of space curves are parametrized by 
the tangent map of the curve and they have rather degenerate singularities \cite{Ishikawa99, Ishikawa12}. 

The purpose of this paper is twofold. First we formulate the normal and tangent maps for 
\lq\lq frontals", a generalized submanifolds, 
or immersions with singularities but with well-defined tangent spaces, giving 
some basics for differential geometric and dynamical study of frontals. 
Second we show a new aspect to provide an inter-relation of \lq\lq normal and tangent". 

In \S \ref{Tangent maps to frontals and their parallels}, 
we define and study locally on normal maps and tangent 
maps of \lq\lq{frontals}". 
Moreover we give a result showing that the 
parallels of the tangent map to a frontal is right equivalent to the tangent map of a frontal under some natural conditions 
(Theorem \ref{parallel-of-tangent}). 
All results of our paper in \S \ref{Tangent maps to frontals and their parallels}
are proved in \S \ref{Proofs of results}. 

In this paper all manifolds and mappings are assumed to be of class $C^\infty$ unless otherwise stated. 

The author would like to thank organizers of DIFF2020 and MOCS2020 for their continuous efforts. 

The author is grateful to the anonymous referee for his/her comment to indicate the errors  
in the previous draft and to show an essential example about it. 
The comment has been very helpful for the author to revise the paper. 

\section{Tangent maps to frontals and their parallels} 
\label{Tangent maps to frontals and their parallels} 

First we introduce the notion of frontals. Let $U$ be any $n$-dimensional manifold. 

\bef
{\rm 
A map-germ $f : (U, a) \to \R^{n+p}$ is 
called a 
{\it frontal} if there exists a smooth ($= C^\infty$) $n$-plane family, 
called a {\it Legendrian lift}, 
$$
(U, a) \ni t \mapsto \widetilde{f}(t) \subset T_{f(t)}\R^{n+p}
$$ 
along $f$ such that $f_*(T_tU) \subseteq \widetilde{f}(t)$. 

A Legendrian lift $\widetilde{f}$ of $f$ is regarded as an integral lift 
$\widetilde{f} : (U, a) \to \Gr(n, T\R^{n+p})$ of $f$ for the canonical distribution $D \subset 
T\Gr(n, T\R^{n+p})$ of the Grassmannian bundle over $\R^{n+p}$: 
For $(x, V) \in \Gr(n, T\R^{n+p})$ with $x \in \R^{n+p}, V \subset T_x\R^{n+p}$, $\dim(V) = n$, 
we set $D_{(x, V)} := \pi_*^{-1}(V)$, where $\pi_* : T_{(x, V)}\Gr(n, T\R^{n+p}) \to T_x\R^{n+p}$ 
is the differential of the projection $\pi : \Gr(n, T\R^{n+p}) \to \R^{n+p}$, $\pi(x, V) := x$. 

A frontal is often written as the pair $(f, \widetilde{f})$, if $\widetilde{f}$ is a fixed Legendrian lift. 

A map $f : U \to \R^{n+p}$ is called a {\it frontal} if the germ of $f$ at any point 
$a \in U$ is frontal. 
}
\enf

Then, given a Legendrian lift $\widetilde{f}$ over $U$, 
the pull-back bundle $f^*(T\R^{n+p})$ is decomposed into the 
sum $T_f \oplus N_f$ of 
the {\it tangent bundle} $T_f$ of rank $n$ 
and the {\it normal bundle} $N_f$ of rank $p$ 
of $f$ over $U$: 
$$
T_{f, t} := \widetilde{f}(t), \ \ 
N_{f, t} := T_{f, t}^\perp, 
\ \ T_{f(t)}\R^{n+p} = T_{f, t} \oplus N_{f, t}. 
$$

The condition that $f$ is a frontal, is equivalent to the existence of local 
orthonormal normal frame $\{ \nu_1, \dots, \nu_p\}$. In other words, 
being a frontal is equivalent to be locally framable. 

\

\bef
\label{proper-def}
{\rm 
A map-germ $f : (U, a) \to \R^{n+p}$ is 
called {\it proper} if 
the singular locus
$$
\Sigma(f) := \{ t \in (U, a) \mid \rank(T_tf : T_tU \to T_{f(t)}\R^{n+p}) < n \} 
$$ 
of $f$ has no interior point nearby $a \in \R^n$. 
}
\enf

\ber
{\rm 
Any Legendrian lift $\widetilde{f}(t)$ of a frontal $f$ necessarily coincides with $T_tf(T_tU)$ 
provided $t \not\in \Sigma(f)$. Therefore we have that 
the Legendrian lift $\widetilde{f}$ of a proper frontal-germ 
is uniquely determined (\cite{Ishikawa18} Proposition 6.2). 
}
\enr

\bep
\label{principal}
{\rm (\cite{Ishikawa20} Lemma 2.3.)} 
Let $f : (\R^n, 0) \to \R^{n+p}$ be a map-germ. 
If $f$ is a frontal-germ, then the Jacobi ideal $J_f$, that is generated by 
$n$-minor determinants of the Jacobian matrix of $f$, is a
principal ideal, i.e. it is generated by one element. 

Conversely, if the Jacobi ideal $J_f$ is a principal ideal and 
$f$ is proper, then $f$ is a frontal-germ. 
\enp

\bef
\label{normal-tangent-def}
{\rm 
The {\it normal map} ${\mathrm{Nor}}(f) : N_f \to \R^{n+p}$  (resp. 
the {\it tangent map} ${\mathrm{Tan}}(f) : T_f \to \R^{n+p}$) to a frontal $f$
is defined by $(t, \nu) \mapsto f(t) + \nu$ \ (resp. $(t, \tau) \mapsto f(t) + \tau$\ ). 
}
\enf

Then it is naturally expected that the singularity and geometry of 
the map ${\mathrm{Nor}}(f)$ (resp. ${\mathrm{Tan}}(f)$) 
reflects those of the frontal $f$. 

\

Let $f : (\R^n, 0) \to \R^{n+p}$ be a frontal-germ and 
${\mbox{\rm Nor}}(f) : N_f \to \R^{n+p}$ the normal map of $f$. 
If $\{ \nu_1, \dots, \nu_p\}$ is a local frame of the normal bundle $N_f$, then 
the normal map is represented by ${\mbox{\rm Nor}}(f) : (\R^n, 0)\times \R^p \to 
\R^{n+p}$, 
$$
{\textstyle 
{\mbox{\rm Nor}}(f)(t, u) = f(t) + \sum_{i=1}^p u_i \nu_i(t), \quad t \in (\R^n, 0), u \in \R^p. 
}
$$
The right equivalence class of ${\mbox{\rm Nor}}(f)$ is independent of the choice of 
$\{ \nu_1, \dots, \nu_p\}$. 

The normal map has a lift 
$$
\widetilde{\mathrm{Nor}}(f) : N_f \to T\R^{n+p} \cong T^*\R^{n+p}, 
$$
which is defined by
$$
\widetilde{\mathrm{Nor}}(f)(t, \nu) = (f(t) + \nu; \ \nu^*), 
$$
where $\nu^* \in T_{f(t) + \nu}^*\R^{n+p}$ is defined by 
$\nu^*(w) = \nu\cdot w, (w \in T_{f(t) + \nu}\R^{n+p})$. 
Then, for the canonical symplectic form $\Omega$ 
on $T\R^{n+p} \cong T^*\R^{n+p}$, via the Euclidean metric, 
the pull-back of $\Omega$ by $\widetilde{\mathrm{Nor}}(f)$ vanishes, i.e. 
the normal map lifts a Lagrangian map. 

\bef
{\rm 
Let $f : (\R^n, 0) \to \R^{n+1}$ be a frontal hypersurface and 
$\nu : (\R^n, 0) \to S^{n}$ a unit normal field along $f$. 
The mappings $f + u\nu : (\R^n, 0) \to \R^{n+1}, (u \in \R)$ 
have the common normal line fields with the original $f$, 
and they are called {\it parallel hypersurfaces} to $f$. 
}
\enf

Let $f : (\R^n, 0) \to \R^{n+p}$ be a frontal, $r > 0$ and consider a $(p-1)$-sphere bundle 
$
(N_f)_r := \{ (t, \nu) \in N_f \mid \Vert \nu \Vert = r \}
$ over $(\R^n, 0)$, $\dim(N_f)_r = n+p-1$. 
Then the restriction of the normal map
$$
{\mbox{Can}}(f) := {\mbox{\rm Nor}}(f)\vert_{(N_f)_r} : (N_f)_r \to \R^{n+p}
$$
to $(N_f)_r$ is called a {\it canal hypersurface}.  It has the Legendrian lift 
$$
{\widetilde{\mbox{Can}}}(f) : (N_f)_r \longrightarrow P(T\R^{n+p}) \cong P(T^*\R^{n+p}), 
$$
which is defined by ${\widetilde{\mbox{Can}}}(f)(t, \nu) := (f(u), t+\nu; [\nu^*])$, 
and the canal hypersurface turns to be a frontal. In fact ${\mbox{Can}}(f)$ 
is regarded as a parallel hypersurface to the frontal hypersurface $F : (N_f)_r \to \R^{n+p}$ 
defined by $F(t, \nu) := f(t)$ and the unit normal $\widetilde{\nu} : (N_f)_r \to S^{n+p-1}$ 
defined by $\widetilde{\nu}(t, \nu) := \frac{1}{r}\nu$. 

\

Now, let us consider the case $p \geq 2$. 
Let $f : (\R^n, 0) \to \R^{n+p}$ be a frontal and $\{ \nu_1, \dots, \nu_p\}$ 
be an orthonormal frame of 
the normal bundle $N_f$. Consider the family of maps 
$f_u(t) := f(t) + \sum_{i=1}^p u_i\nu_1(t)$ with parameter $u$. Then 
all $\nu_1, \dots, \nu_p$ are normal to all $f_u$ if and only if 
$$
\frac{\pa}{\pa t_j}\{ f(t) + \sum_{i=1}^p u_i\nu_i(t)\} \cdot \nu_k = \sum_{i=1}^p u_i\frac{\pa}{\pa t_j}\nu_i(t) \cdot \nu_k = 0, 
$$
for any $1 \leq k \leq p, 1 \leq j \leq n$. Here $\cdot$ means the inner product on $\R^{n+p}$. 
Then we obtain a natural sufficient condition that 
$\frac{\pa}{\pa t_j}nu_i \in T_f \ (1 \leq i \leq p, 
1 \leq j \leq n)$ for that $\{ \nu_1, \dots, \nu_p\}$ is a frame of $N_{f_u}$ for any $u$. 
From this motivation, we formulate as follows: 

\

Let $f : U^n \to \R^{n+p}$ be a frontal and 
$f^*(T\R^{n+p}) = T_f \oplus N_f$ the decomposition into the tangent bundle and 
the normal bundle to $f$ associated to a Legendrian lift $\widetilde{f}$ of $f$. 
For any vector field $\eta$ over $U$ and any vector field $X : U \to T\R^{n+p}$ along $f$, 
we denote by $\nabla_\eta^fX$ the covariant derivative of $X$ by $\eta$ 
induced from the Euclidean metric on $\R^{n+p}$.  
If $\eta(t) = \sum_{i=1}^n a_i(t)\frac{\pa}{\pa  t_i}$ and 
$X(t) = \sum_{j=1}^{n+p} X_j(t) (\frac{\pa}{\pa x_j}\circ f)$, then 
$$
{
\nabla_\eta^fX 
= \sum_{j=1}^{n+p} \sum_{i=1}^n \left\{ a_i\frac{\pa X_j}{\pa t_i}\left(\frac{\pa}{\pa x_j}\circ f\right)+ 
a_iX_j\left(\nabla_{\frac{\pa}{\pa t_i}}^f\frac{\pa}{\pa x_j}\right) \right\}
= \sum_{j=1}^{n+p} \sum_{i=1}^n a_i\frac{\pa X_j}{\pa t_i}\left(\frac{\pa}{\pa x_j}\circ f\right), 
}
$$
since $\nabla_{\frac{\pa}{\pa t_i}}^f\frac{\pa}{\pa x_j} = 0$ in the Euclidean case. 
Note that $\nabla_\eta^fX$ is a section of $f^*(T\R^{n+p})$, i.e. a vector field 
along $f$. 

By the decomposition $f^*(T\R^{n+p}) = T_f \oplus N_f$, we write 
$$
\nabla_\eta^fX = \nabla^\top_\eta X + \nabla^\perp_\eta X. 
$$
Then $\nabla^\perp$ (resp. $\nabla^\top$) defines a connection on the 
vector bundle $N_f$ (resp. $T_f$). 

\bef
{\rm 
The induced connection $\nabla^\perp$ on $N_f$ is called 
the {\it normal connection} or {\it Van der Waerden-Bortolotti connection} of the framed frontal $f$. 
The connection $\nabla^\top$ on $T_f$ is called 
the {\it tangential connection} or {\it Levi-Civita connection} of $f$. 

We call the frontal $f$ {\it normally flat} (resp. {\it tangentially flat}) 
if the normal connection on $N_f$ (resp. the tangential connection on $T_f$)
is flat, i.e. 
there exists an orthonormal frame $\{ \nu_1, \dots, \nu_p\}$ of $N_f$, 
(resp. $\{ \tau_1, \dots, \tau_n\}$ of $T_f$) such that 
$\nabla^\perp_\eta\nu_i = 0, 1 \leq i \leq p$ 
(resp. $\nabla^\top_\eta\tau_j = 0, 1 \leq j \leq n$), 
for any vector field $\eta$ on $U$. 

We call $\{ \nu_1, \dots, \nu_p\}$ 
a {\it normally parallel orthonormal frame} of $N_f$ 
or briefly a {\it Bishop frame} of the normally flat frontal (see \cite{Bishop}). 
We call $\{ \nu_1, \dots, \nu_p\}$ a {\it tangentially parallel orthonormal frame} of $T_f$ 
of the tangentially flat frontal. 
}
\enf

\

Then we will show the followings. 

\bel
\label{hypersurface-normally-flat}
{\rm 
Any frontal hypersurface $f : (\R^n, 0) \to \R^{n+1}$ is normally flat. In fact, any unit normal is 
normally parallel. 
}
\enl

\bel
\label{curve-normally-flat}
{\rm 
Any frontal curve $f : (\R, 0) \to \R^{1+p}$ is normally flat and tangentially flat. 
}
\enl

Let $f : (\R, 0) \to \R^{1 + p}$ be a frontal curve. 
An orthonormal frame $\nu_1, \dots, \nu_p$ of $N_f$ is 
a Bishop frame 
if and only if $(\nu_j)_u \in T_\gamma, (1 \leq j \leq p)$ (\cite{Bishop}). 
Then the curves 
$f + \sum_{j=1}^p u_j\nu_j$ are called {\it parallel curves} to $f$. 
They are all frontals and common Bishop frame with $f$. 

Parallels are basic and interesting objects to be studied in the cases of both hypersurfaces and curves 
(\cite{Bruce, FH, HT, LPTY}).

\ber
{\rm 
Let $f : (\R, 0) \to \R^{1+p}$ be a frontal curve. 
Let $\tau(t)$ be a unit section of $T_f$ and $\nu_1(t), \dots, \nu_p(t)$ a 
normally parallel orthonormal frame of $N_f$. 
Then we have the structure equation of the framed curve $(f; \tau, \nu_1, \dots, \nu_p)$ 
given by $f'(t) = a(t) \tau(t)$ and 
$$
\left(
\begin{array}{cccc}
\tau'(t), \ \nu_1'(t), \ \dots, \nu_p'(t)
\end{array}
\right)
= 
\left(
\begin{array}{cccc}
\tau(t), \ \nu_1(t), \ \dots, \ \nu_p(t)
\end{array}
\right)
\left(
\begin{array}{cccc}
0 & -\kappa_1(t) & \cdots & -\kappa_p(t)
\\
\kappa_1(t) & 0 & \cdots & 0
\\
\vdots & \vdots & \ddots & \vdots 
\\
\kappa_p(t) & 0 & \cdots & 0
\end{array}
\right)
$$
The functions $a, \kappa_1, \dots, \kappa_p$ form a complete system of invariants 
of frontal curves up to congruence (cf. \cite{FT}). 
}
\enr

\ber
{\rm 
For smooth surfaces in $\R^4$, the normal flatness is characterized by the vanishing of 
the normal curvature (see \cite{Little} for instance). 
Note that the normal curvature plays important role in the study of submanifolds from the 
view points of singularity theory (see \cite{MFR, IFRT}). The notion of normal curvature is naturally generalized to frontals: 
Let us consider, for instance, a frontal surface $f : (\R^2, 0) \to \R^4$ in $\R^4$ with 
the decomposition $T_f \oplus N_f$ of $f^*(T\R^4)$. Let $e_1, e_2, e_3, e_4$ be an orthonormal frame of 
of $f^*(T\R^4)$ such that $e_1, e_2$ (resp. $e_3, e_4$) generate $T_f$ (resp. $N_f$). 
Let $\{ \omega_1, \omega_2, \omega_3, \omega_4 \}$ be the dual frame of 
$\{ e_1, e_2, e_3, e_4\}$. 
The normal connection $\nabla^\perp$ determines the connection $1$-form 
$\omega_{34}$ by the equality $\nabla^{\perp}_{\eta}e_3 = \omega_{34}(\eta)e_4$ for any 
vector field $\eta$ over $(\R^2, 0)$. Then the curvature $2$-form $\Omega_{34}$ is defined 
by $\Omega_{34} = d\omega_{34}$ and the normal curvature function $N$ over $(\R^2, 0)$ 
of $f$ is defined by $\Omega_{34} = - K \omega_1\wedge\omega_2$. 
Then the frontal is normally flat if and only if $K = 0$. 

Normally flat submanifolds of arbitrary codimensions and dimensions are studied as 
an important class in the differential geometry of submanifolds (see \cite{Terng, BCO}). 
In this paper we study, as one of possible first steps, the local geometry of normally flat frontals 
via tangent and normal maps. 
}
\enr

\bef
{\rm 
Let $f : (\R^n, 0) \to \R^{n+p}$ be a normally flat frontal, and 
$\{ \nu_1, \dots, \nu_p\}$ a Bishop frame of $N_f$. Then we define 
the {\it parallels} to normally flat frontal $f$ by 
$$
{\textstyle 
{\mbox{\rm Pal}}_\nu(f) := 
f(t) + \nu(t), \quad \nu(t) := \sum_{i=1}^p u_i \nu_i(u), \quad (u \in \R^p). }
$$
Then the parallels are normally flat and have the same Bishop frame with $f$. 
}
\enf

Next we turn to tangent maps. 
Let $f : (\R^n, 0) \to \R^{n+p}$ be a frontal and 
${\mbox{\rm Tan}}(f) : T_f \to \R^{n+p}$ the tangent map of $f$ 
(Definition \ref{normal-tangent-def}). 
If $\{ \tau_1, \dots, \tau_n\}$ is a local frame of the tangent bundle $T_f$, then 
the tangent map is represented by ${\mbox{\rm Tan}}(f) : (\R^n, 0)\times \R^n \to 
\R^{n+p}$, 
$$
{\textstyle 
{\mbox{\rm Tan}}(f)(t, s) = f(t) + \sum_{i=1}^n s_i \tau_i(t), \quad t \in (\R^n, 0), s \in \R^n. 
}
$$
The right equivalence class of ${\mbox{\rm Tan}}(f)$ is independent of the choice of 
$\{ \tau_1, \dots, \tau_n \}$.

\

Recall that a map-germ is called proper if 
its singular locus is nowhere dense (Definition \ref{proper-def}). 

\bet
\label{tangent-nomally-flat}
{\rm (Tangent maps of normally flat frontals are normally flat.)} \ Let $n \leq p$ and 
$f : (\R^n, 0) \to \R^{n+p}$ be a normally flat frontal. 
Suppose $\Tan(f) : (\R^{n}, 0)\times\R^n \to \R^{n+p}$ is a proper frontal of $f$, 
and the restriction $N_{\Tan(f)}\vert_{\R^n\times\{0\}}$ is flat. 
Then $\Tan(f)$ is a normally flat frontal. 
\ent

The conditions on $\Tan(f)$ and $N_{\Tan(f)}\vert_{\R^n\times\{0\}}$ 
are required to guarantee that 
the Legendrian lift of $\Tan(f)$ is well-controlled in the proof. 
The author does not know whether these conditions can be weakened. 
Note that, if $n = 1$, then 
$N_{\Tan(f)}\vert_{\R^n\times\{0\}}$ is necessarily flat (Lemma \ref{flatness-of-connection}).

\bee
{\rm
Even if $f$ is normally flat frontal, $\Tan(f)$ is not necessarily a frontal. 
For example, let $f : (\R, 0) \to \R^3$ be a $C^\infty$ curve defined by 
$$
f(t) = (x(t), y(t), z(t)) = \left\{ 
\begin{array}{c}
 (\exp(-1/t^2), 0, -t^2) \ (t \geq 0), 
\\
(0, \exp(-1/t^2), -t^2) \ (t \leq 0). 
\end{array}\right. 
$$
Then $f$ is a frontal curve and $\Tan(f) : (\R^2, 0) \to \R^3$ is not a frontal. 
\begin{center} 
\includegraphics[width=4truecm, height=4truecm, clip, 
bb= 149 536 440 798]{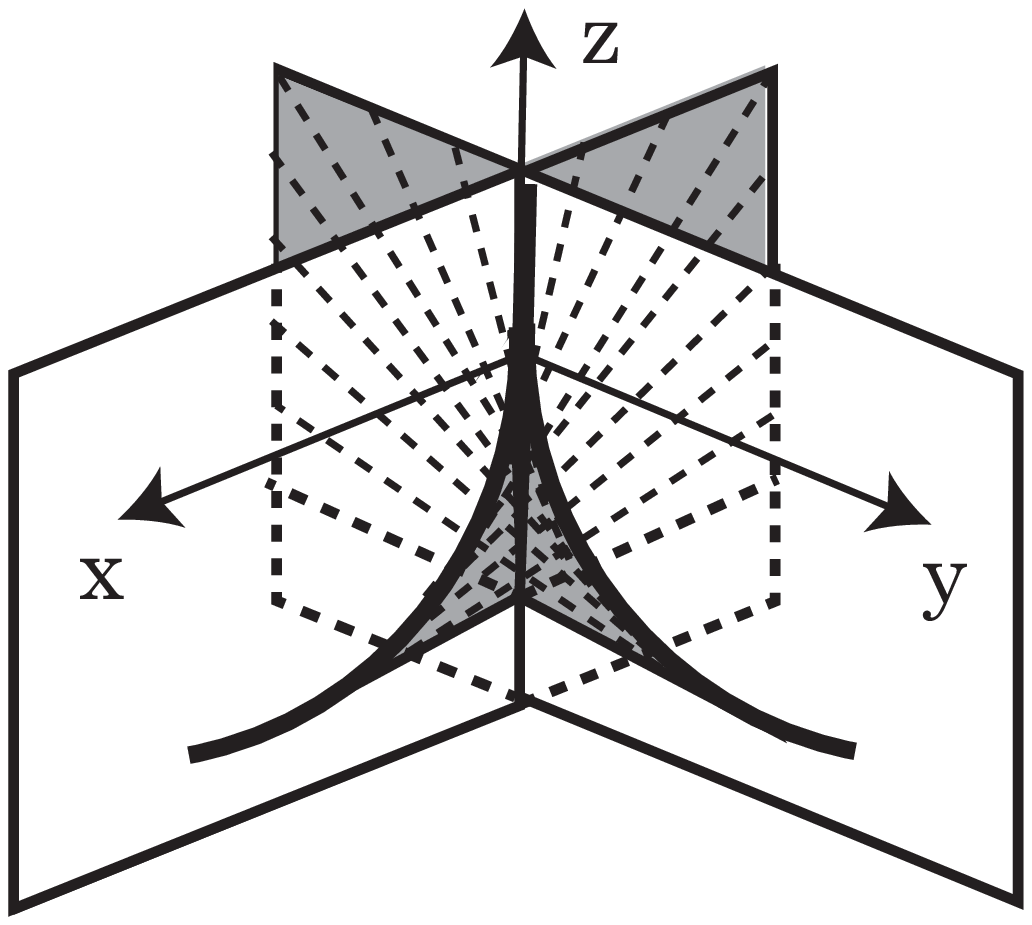}  
\\
A frontal curve with the non-frontal tangent surface.
\end{center}
}
\ene

We give here a sufficient condition on a curve to be a frontal. 
To show it we define, for a curve $f : (\R^1, 0) \to \R^{1+p}$, the $(n+p)\times k$-Wronskian matrix of $f$ by 
$$
W_k(f)(t) := 
\left( \dfrac{df}{dt}(t), \dfrac{d^2f}{dt^2}(t), \dfrac{d^3f}{dt^3}(t), \cdots, \dfrac{d^kf}{dt^k}(t)\right), 
$$
for any positive integer $k$. Note that the rank of $W_k(f)(t)$ is independent of the choice of 
local coordinates of $(\R^1, 0)$. 

\bel
\label{sufficient condition}
Let $f : (\R, 0) \to \R^{1+p}$ be a curve. 
Suppose $W_k(f)(0)$ at $0$ is of rank $\geq 2$ for a sufficiently large $k$. 
Then $f$ is a frontal and its tangent map $\Tan(f) : (\R, 0)\times \R \to \R^{1+p}$ is frontal. 
\enl

In this paper we show 
an inter-relation of \lq\lq normal and tangent" of a frontal, 
as the main theorem of this paper. To state it, we prepare some geometric notions on 
frontal curves ant their tangent maps. 

\bef
{\rm 
Let $f : I \to \R^{1+p}$ be a frontal curve defined on an interval $I$. 
A point $t_0 \in I$ is called an {\it inflection point} 
of $f$ if there exists a unit frame $\tau : I \to S^p$ of $T_f$ such that 
$t_0$ is a singular point of $\tau$, regarded as a curve in $S^p$, 
i.e. $\tau'(t_0) = 0$. 
}
\enf

Suppose $\Tan(f) : I\times \R \to \R^{1+p}$ is a proper frontal and 
$T_{\Tan(f)}$ is the tangent bundle of $\Tan(f)$. Let $\tau$ be a unit section 
of $T_f$ regarded as a unit section of $\Tan(f)$ along $I \times\{ 0\}$. Note that 
$\Tan(f)(t, 0) = f(t)$. 
Take an orthonormal frame $\{ \tau, \mu\}$ of $T_{\Tan(f)}$ along $I\times\{ 0\}$. 
Then there exists unique function $\kappa : I \to \R$ with $\tau'(t) = \kappa(t)\mu(t)$, 
which may be called the {\it curvature function} of the frontal curve $f$ with respect 
the frame $\{ \tau, \mu\}$. Remark that $t_0 \in I$ is an inflection point of $f$ if 
and only if $\kappa(t_0) = 0$. 
Take a normally parallel orthonormal frame $\{ \nu_1, \dots, \nu_{p-1}\}$ of 
$N_{\Tan(f)}$ along $I\times\{ 0\}$. Then there exist functions 
$\ell_1, \dots, \ell_{p-1} : I \to \R$ uniquely, which may be called {\it torsion functions} of $f$, 
such that $\mu'(t) = - \kappa(t)\tau(t) + 
\sum_{i=1}^{p-1} \ell_i(t)\nu_i(t)$. 
The frame $\{ \nu_1, \dots, \nu_{p-1}\}$ extends to a 
normally parallel orthonormal frame of $N_{\Tan(f)}$ over $I\times\R$ 
as $\nu_i(t, s) = \nu_i(t)$, by the constancy of $N_{\Tan(f)}$ along $t\times \R$. 

Now let us consider parallels ${\mbox{\rm Pal}}_\nu(\Tan(f))$ of $\Tan(f)$ 
by normals $\nu = \sum_{i=1}^{p-1} u_i\nu_i$.

\bet
\label{parallel-of-tangent}
{\rm  (Parallels to tangent maps of frontal curves are tangent maps of frontal curves.)} \ 
Let $p \geq 2$. 
Let $f : (\R, 0) \to \R^{1+p}$ be a frontal curve-germ without inflection point, 
and suppose $\Tan(f) : (\R, 0)\times \R \to \R^{1+p}$ is a proper frontal. 
Then any parallel ${\mbox{\rm Pal}}_\nu(\Tan(f))$, $\nu = \sum_{i=1}^{p-1} u_i\nu_i$, 
to $\Tan(f)$ is right equivalent to the tangent map $\Tan(g)$ 
for a frontal curve, called the directrix or the edge of regression, $g : (\R, 0) \to \R^{1+p}$. 
In fact $g$ is given by 
$$
g(t) = f(t) + \sum_{i=1}^{p-1} u_i\left(\frac{\ell_i(t)}{\kappa(t)}\,\tau(t) + \nu_i(t)\right). 
$$
\ent

\ber
{\rm 
In the first draft of the paper, 
the author claimed erroneously that parallels to tangent maps of frontals are tangent maps of parallels to frontals. 
The above Theorem \ref{parallel-of-tangent} provides a corrected and weakened revision and 
the curve $g$ is not necessarily a parallel to $f$ as in the following example, 
which is given by the referee to the author. 
}
\enr

\bee
{\rm 
Let $f : (\R, 0) \to \R^3$ be defined by $f(t) = (t, \frac{t^2}{2}, \frac{t^3}{6})$.  
Then $\Tan(f) : (\R, 0)\times \R \to \R^3$ 
is given by $\Tan(f)(t, s) = (t+s, \frac{t^2}{2} + st, \frac{t^3}{6} + s\frac{t^2}{2})$. 
The unit frame $\tau$ of $T_f$ is given by $\tau(t) = \frac{2}{2+t^2}(1, t, \frac{t^2}{2})$. 
The orthonormal frame of $T_{\Tan(f)}$ is given by $\{ \tau, \mu\}$ with 
$\tau(t, s) = \tau(t)$ and $\mu(t, s) = \frac{1}{2+t^2}(-2t, 2-t^2, 2t)$. 
The vector field $\nu(t, s) = \nu_1(t, s) := \frac{1}{2+t^2}(t^2, -2t, 2)$ is a unit normal 
to $\Tan(f)$. 
The parallel to $\Tan(f)$ along $\nu$ is defined by 
${\mbox{\rm Pal}}_{u\nu}(\Tan(f))(t, s) = \Tan(f)(t, s) + u\nu(t, s) 
= (t+s+ \frac{ut^2}{2+t^2}, \frac{t^2}{2} + st - \frac{2ut}{2+t^2}, \frac{t^3}{6} + s\frac{t^2}{2} + \frac{2u}{2+t^2})$, 
$u \in \R$. 
In this example, both curvature $\kappa$ and torsion $\ell = \ell_1$ of $f$ 
are calculated as $\kappa(t) = \ell(t) = \frac{2}{2+t^2}$ 
and thus $\frac{\ell(t)}{\kappa(t)}\tau + \nu(t) = (u, 0, u)$. Then 
${\mbox{\rm Pal}}_{u\nu}(\Tan(f))$ is right equivalent to $\Tan(g)$, 
where the directrix $g$ is obtained by $g(t) = (t + u, \frac{t^2}{2}, \frac{t^3}{6} + u)$. 
Moreover we see that $g$ is a parallel to $f$ if and only if $u = 0$. 
}
\ene

\bee
{\rm In Theorem \ref{parallel-of-tangent}, the condition that $f$ has no inflection point is necessary: 
Let $f : (\R, 0) \to \R^3$ be defined by $f(t) = (t, \frac{t^3}{6}, \frac{t^4}{24})$. 
Then $f$ has an inflection point at $t_0 = 0$. 
We have $\Tan(f)(t, s) = (t+s, \frac{t^3}{6} + s\frac{t^2}{2}, \frac{t^4}{24} + s\frac{t^3}{6})$  
and the unit normal to $\Tan(f)$ is given by 
$\nu(t, s) = (\nu_1(t), \nu_2(t), \nu_3(t)) = \frac{12}{\sqrt{t^6+36t^2+144}}(\frac{t^3}{12}, -\frac{t}{2}, 1)$. 
The singular locus of the parallel ${\mbox{\rm Pal}}_{u\nu}(\Tan(f))$ is 
given by $st - u(\frac{t^2}{2}\nu_1'(t) - \nu_2'(t)) = 0$. Note that $\nu_2'(0) = - \frac{1}{2} \not= 0$. 
The parallel is a tangent surface outside of $t = 0$. However 
the directrix diverges to infinity as $t$ tends to $0$, provided that $u \not= 0$. 
}
\ene

The generic classification of singularities in parallels of the tangent surfaces to 
space frontal curves will be given in a forthcoming paper. 

\section{Proofs of results}
\label{Proofs of results}

First we remark

\bel
Let $f : (\R^n, 0) \to \R^{n+p}$ 
be a normally flat frontal and $\nu : (\R^n, 0) \to T\R^{n+p}$ 
a normally parallel vector field over $f$. 
Then $\nu$ is uniquely determined by the value $\nu(0)$ at $0$. 
The vector space of normally parallel vector fields along $f$ is 
isomorphic to $\R^p$. 
\enl

\Proof
Let $\{\nu_1, \dots, \nu_p\}$ be a normally parallel frame. 
Write $\nu = \sum_{i=1}^p a_i \nu_i$ for some function-germs $a_i$ on $(\R^n, 0)$. 
Take any vector field $\eta$ over $(\R^n, 0)$. Then 
$\nabla_\eta^f \nu = \sum_{i=1}^p (\eta a_i)\nu_i + \sum_{i=1}^p a_i \nabla_\eta^f \nu_i$. 
Since $\nabla_\eta^f \nu_i, 1 \leq i \leq p$ are tangential, 
we have 
$\nabla_\eta^\perp \nu= \sum_{i=1}^p (\eta a_i)\nu_i$. 
Hence $\nu$ is normally parallel if and only if $a_i (1 \leq i \leq p)$ are all constant. 
Then $\nu$ is uniquely determined by $\nu(0) = \sum_{i=1}^p a_i\, \nu_i(0)$. 
\qed

\

Next we recall 

\bel
The normal connection $\nabla^\perp$ on $N_f$  and the tangential connection $\nabla^\top$
are metric preserving: 
$$
\begin{array}{ccc}
\xi(\nu \cdot \mu) & = & (\nabla^\perp_\xi \nu)\cdot \mu + \nu \cdot(\nabla^\perp_\xi \mu), 
\vspace{0.1truecm}
\\
\xi(\tau \cdot \rho) & = & (\nabla^\top_\xi \tau)\cdot \rho + \tau \cdot(\nabla^\top_\xi \rho), 
\end{array}
$$
for any $\xi \in \Gamma(TU), \nu, \mu \in \Gamma(N_f), \tau, \rho \in \Gamma(T_f)$. Here $\cdot$ means the inner product on $\R^{n+p}$. 
\enl

\Proof
The first equality is obtained as follows: 
$$
\begin{array}{ccl}
\xi(\nu \cdot \mu) & = & \nabla_\xi(\nu \cdot \mu) = (\nabla_\xi\nu)\cdot\mu + \nu\cdot(\nabla_\xi\mu) 
\\
& = & (\nabla^\perp_\xi \nu + \nabla^\top_\xi \nu)\cdot \mu 
+ \nu \cdot(\nabla^\perp_\xi \mu + \nabla^\top_\xi \nu) 
\\
& = & (\nabla^\perp_\xi \nu)\cdot \mu + \nu \cdot(\nabla^\perp_\xi \mu). 
\end{array}
$$
The second equation is obtained similarly. 
\qed

\

Lemma \ref{hypersurface-normally-flat} and 
Lemma \ref{curve-normally-flat} follow from the following general well-known result: 

\bel
\label{flatness-of-connection}
Let $U$ be an $n$-dimensional manifold. Let $(E, g)$ be a 
metric vector bundle of rank $p$ on $U$ and $\nabla$ a metric-preserving connection 
on $E$: $\nabla : \Gamma(E) \to \Gamma(T^*U\otimes E)$. 
If $n = 1$ or $p = 1$, then $\nabla$ is flat, i.e., 
locally there exists $\nabla$-parallel $g$-orthogonal frame of $E$. 
\enl

\Proof 
A connection if flat if and only if its curvature form vanishes identically (see \cite{KN} Theorem 9.1 of chapter II, for instance). The curvature form of a metric-preserving connection is, locally, 
a skew-symmetric $p\times p$-matrix of $2$-forms of $n$-variables, and therefore it vanishes 
identically if $n = 1$ or $p = 1$. 
\qed

\

Here we give, to make sure, alternative and direct proofs of Lemma \ref{hypersurface-normally-flat} and 
Lemma \ref{curve-normally-flat}. 

\

\noindent
{\it Proof of Lemma \ref{hypersurface-normally-flat}.} \
Let $f : (\R^n, 0) \to \R^{n+1}$ be a frontal hypersurface and 
$\nu : (\R^n, 0) \to T\R^{n+1}$ any unit normal field along $f$. 
Then $\nu\cdot\nu = 1$ implies $\nu_{u_i}\cdot\nu = 0, 1 \leq i \leq n$. 
This means that each $\nu_{u_i}, (1 \leq i \leq n)$ is tangential. 
Therefore $\nu$ is normally parallel,  $\nabla^\perp\nu = 0$, and $f$ is normally flat. 
\qed

\

\noindent
{\it Proof of Lemma \ref{curve-normally-flat}.} \
Let $f : (\R, 0) \to \R^{1+p}$ be a frontal curve. 

First let us show that $f$ is tangentially flat. Let $\tau : (\R, 0) \to T\R^{1+p}$ 
the unit frame of $T_f$. Then $\tau\cdot\tau = 1$ 
implies $\tau'\cdot\tau = 0$. 
This means that the differential $\tau'$ is normal. 
Therefore $\tau$ is tangentially parallel,  $\nabla^\top\tau = 0$, and $f$ is tangentially flat. 

Next we show that $f$ is normally flat. 
Let $\widetilde{f}$ be a Legendrian lift of the frontal $f$. 
We may suppose $\widetilde{f}(0)$ is generated by the vector ${ }^t(1, 0, \dots, 0)$ 
for a system of Euclidean coordinates. 
Let ${ }^t(1, a_1(t), \dots, a_p(t))$ be a generator of $\widetilde{f}(t)$ with 
$a_j(0) = 0, 1 \leq j \leq p$. 
Note that in this case a unit tangent vector $\tau$ is obtained by 
$$
\tau(t) = \frac{1}{1+\sum_{i=1}^p a_i^2}{ }^t(1, a_1(t), \dots, a_p(t)). 
$$

Let $\nu : (\R, 0) \to T\R^{1+p}$ any normal field along $f$. We seek the condition that $\nu$ is 
normally parallel. 
Set 
$\nu(t) = { }^t(b_0(t), b_1(t), \dots, b_p(t))$. Then 
$b_0 + b_1a_1 + \cdots + b_pa_p = 0$. Therefore 
$$
\nu = b_1\,{ }^t(-a_1, 1, 0, \dots, 0, 0) + b_2\,{ }^t(-a_2, 0, 1, \dots, 0, 0)+ \cdots + b_p\,{ }^t(-a_p, 0, 0, \dots, 0, 1). 
$$
The normal field $\nu$ is normally parallel if and only if there exists a function $\lambda = \lambda(t)$ 
such that 
$\nu_t = \lambda\, { }^t(1, a_1, \dots, a_p)$. 
The quality is equivalent to that 
$$
(b_j)_t = \lambda a_j \ (1 \leq j \leq p), \quad \lambda = - \sum_{j=1}^p (b_j)_t a_j - 
\sum_{j=1}^p b_j (a_j)_t. 
$$
Then, by a straightforward computation, 
the condition that $\nu$ is normally parallel is given by the linear 
ordinary differential equation 
{\small 
$$
\left(
\begin{array}{c}
(b_1)_t
\\
(b_2)_t 
\\
\vdots
\\
(b_p)_t 
\end{array}
\right)
= 
\left(
\begin{array}{cccc}
1 + a_1^2 & a_1a_2 & \cdots & a_1a_p
\\
a_2a_1 & 1 + a_2^2 & \cdots & a_2a_p
\\
\vdots & \vdots & \ddots & \vdots 
\\
a_pa_1 & a_pa_2 & \cdots & 1 + a_p^2
\end{array}
\right)^{-1}
\left(
\begin{array}{cccc}
a_1(a_1)_t & a_1(a_2)_t & \cdots & a_1(a_p)_t
\\
a_2(a_1)_t & a_2(a_2)_t & \cdots & a_2(a_p)_t
\\
\vdots & \vdots & \ddots & \vdots 
\\
a_p(a_1)_t & a_p(a_2)_t & \cdots & a_p(a_p)_t
\end{array}
\right)
\left(
\begin{array}{c}
b_1
\\
b_2
\\
\vdots
\\
b_p 
\end{array}
\right)
$$
}
on $b_1(t), \dots, b_p(t)$, with 
$b_0 = - b_1a_1 - \cdots - b_pa_p$. 

Then for any given initial value 
${ }^t(0, \mu_{10}, \dots, \mu_{p0}) \in (N_f)_0 \subset T_{f(0)}\R^{1+p}$, 
there exists unique normally parallel section $\nu$ of $N_f$ 
such that $\nu(0) = { }^t(0, \mu_{10}, \dots, \mu_{p0})$. 
Since $\nabla^\perp$ preserves the metric, given any orthonormal basis 
of $(N_f)_0$, 
there exists unique normally parallel frame, i.e. Bishop frame,  of $N_f$ with the initial value. 
In particular, for each $i (1 \leq i \leq p)$, let $\nu_i$ be the normally parallel section of $N_f$ for the initial value $\mu_{j0} = \delta_{ij} (1 \leq j \leq p)$. 
Then $\nu_1, \dots, \nu_p$ form a normally parallel orthonormal frame of $N_f$. 
\qed

\

To show Theorem \ref{tangent-nomally-flat}, we prepare 

\bel
\label{constancy-of-tangent}
Let $f : (\R^n, 0) \to \R^{n+p}$ be a frontal with $n \leq p$. 
Suppose $F = \Tan(f) : (\R^n, 0)\times\R^n \to \R^{n+p}$ 
is a proper frontal with unique Legendrian lift $\widetilde{F} : 
(\R^n, 0)\times\R^n \to \Gr(n, T\R^{n+p})$. 
Then, we have 

{\rm (1)} $T_f(t, 0) \subset T_F(t, 0)$ for $(t, 0) \in (\R^n, 0)\times\{ 0\}$. 

{\rm (2)} 
For each $(t, s) \in (\R^n, 0)\times\R^n, s \not= 0$, 
the plane field $\widetilde{F}$ is constant along the line $(t, \lambda s), \lambda \in \R$ 
with the direction, 
via the parallel translation in the affine space $\R^{n+p}$ along the line, in other words, 
$(\mbox{\rm{tr}}_{(0, s)})_*\widetilde{F}(t, 0) = \widetilde{F}(t, s)$ for any $(t, s) 
\in (\R^\times\R^n, (0, 0))$, 
where $\mbox{\rm{tr}}_{(0, s)} : \R^{n+p} \to \R^{n+p}$, $\mbox{\rm{tr}}_{(0, t)}(t, s') := (t, s'+s)$. 
\enl

\Proof 
Let $\tau_1(t), \dots, \tau_n(t)$ be a frame of $\widetilde{f}(t)$. 
Set $F = \Tan(f)$, which is defined by $F(t, s) = f(t) + \sum_{i=1}^n s_i\tau_i(t)$. 
Since $F$ is a proper frontal, there exists unique Legendrian lift $\widetilde{F}$, and 
the associated decomposition $F^*(T\R^{n+p}) = T_F \oplus N_F$. 
The image of differential map $F_*(T_{(t, s)}\R^{2n})$ is 
generated by 
$\frac{\pa F}{\pa t_1}, \dots, \frac{\pa F}{\pa t_n}, \frac{\pa F}{\pa s_1}, \dots, \frac{\pa F}{\pa s_n}$, 
Since $\frac{\pa F}{\pa t_j} = \frac{\pa f}{\pa t_j} + \sum_{i=1}^n s_i\frac{\pa \tau_i}{\pa t_j}(u)$, 
$\frac{\pa F}{\pa s_k} = \tau_k$, and $\frac{\pa f}{\pa t_j} \in \left\langle \tau_1, \dots, \tau_n\right\rangle_{{\mathcal E}_{(\R^n, 0)}}$, we have that $F_*(T_{(t, s)}\R^{2n})$ is generated by 
$$
\sum_{i=1}^n s_i\frac{\pa \tau_i}{\pa t_1}(t), \ \dots, \ \sum_{i=1}^n s_i\frac{\pa \tau_i}{\pa t_n}(t), \ 
\tau_1(t), \ \dots, \ \tau_n(t). 
$$
By the integrality condition $F_*(T_{(t, s)}\R^{2n}) \subset \widetilde{F}(t, s)$, we have (1). 

On the other hand we have that the singular locus 
$\Sigma(F)$ of $F$ is invariant of $\R^\times$-action 
$(t, s) \mapsto (t, \lambda s), (\lambda \in \R^\times = \R \setminus \{ 0\}$ and that $\Sigma(F)$ contains the zero-section $(\R^n, 0)\times \{ 0\} 
\subset (\R^n, 0)\times\R^n$. 
Since $F$ is a proper frontal, 
for almost all $(t, s) \in (\R^n, 0) \times(\R^n\setminus \{ 0\})$, we have 
$$
\begin{array}{rcl}
\widetilde{F}(t, \lambda s) & = & 
\left\langle \sum_{i=1}^n \lambda s_i\frac{\pa \tau_i}{\pa t_1}(t), \dots, \sum_{i=1}^n \lambda s_i\frac{\pa \tau_i}{\pa t_n}(t), 
\tau_1(t), \dots, \tau_n(t) \right\rangle_{\R}
\\
& = & 
\left\langle \sum_{i=1}^n s_i\frac{\pa \tau_i}{\pa t_1}(t), \dots, \sum_{i=1}^n s_i\frac{\pa \tau_i}{\pa t_n}(t), 
\tau_1(t), \dots, \tau_n(t) \right\rangle_{\R}, 
\end{array}
$$
for any $\lambda  \in \R \setminus \{ 0\}$. 
Thus we see $\widetilde{F}(t, \lambda s)$ is independent 
of $\lambda$, and by the continuity of 
$\widetilde{F}$ is constant on the line $\{ (t, \lambda s) \mid \lambda \in \R\}$ 
for almost all $(t, s) \in (\R^n, 0) \times(\R^n\setminus \{ 0\})$. 
By the continuity of $\widetilde{F}$ again, 
we have that $\widetilde{F}$ on the line 
$\{ (t, \lambda s) \mid \lambda \in \R\}$, 
for all $(t, s) \in (\R^n, 0) \times(\R^n\setminus \{ 0\})$. 
\qed

\

\noindent
{\it Proof 
of Theorem \ref{tangent-nomally-flat}.} \ 
Let $f : (\R^n, 0) \to \R^{n+p}$ be a normally flat frontal. 
Take an orthonormal parallel frame $\{ \mu_1(t, 0), \dots, \mu_{p-n}(t, 0)\}$ of $N_F$ over $\R^n\times\{ 0\}$. 
Then, by Lemma \ref{constancy-of-tangent} (2), 
$\{ \mu_1, \dots, \mu_{p-n}\}$ is regarded an orthonormal frame of $N_F$ over $\R^n\times\R^n$. 
Then $\nabla_{\pa/\pa t_j} \mu_i(t, 0) \in T_f(t, 0) \subset T_F(t, 0)$ by 
Lemma \ref{constancy-of-tangent} (1). 
Therefore, by Lemma \ref{constancy-of-tangent} (2) again, 
$\nabla_{\pa/\pa t_j} \mu_i(t, s)$ belongs to $T_F(t, s)$. Hence 
$\nabla_{\pa/\pa t_j}^\perp \mu_i = 0, (1 \leq i, j \leq n)$. 
Since $\mu_i$ is independent of $s$, we have $\nabla_{\pa/\pa s_k} \mu_i = 0$ and therefore 
$\nabla_{\pa/\pa s_k}^\perp \mu_i = 0, 1 \leq s \leq p$. Thus we have 
$\{ \mu_1, \dots, \mu_{p-n}\}$ is a Bishop frame along $F$. 
This shows $F = \Tan(f)$ is normally flat. \qed

\

\noindent
{\it Proof of Lemma \ref{sufficient condition}.} 
Let $a_1 = \min\{ k \mid \rank(W_k(f)(0)) = 1\}$ and $a_2 = \min\{ k \mid \rank(W_k(f)(0)) = 2\}$. 
Note that $1 \leq a_1 < a_2$. 
Then there exist a $C^\infty$ coordinate $t$ of $(\R, 0)$ and an affine coordinates 
$x_1, x_2, \dots, x_{1+p}$ such that $f$ is given by 
$$
(x_1\circ f)(t) = t^{a_1}, \ \ (x_2\circ f)(t) = t^{a_2} + o(t^{a_2}), \ \ (x_i\circ f)(t) = o(t^{a_2}) \ (i = 3, \dots, 1+p). 
$$
Then we have a tangential field $\tau(t)$ for $t \not= 0$ by setting $\tau(t) = 
\frac{1}{a_1t^{a_1-1}}f'(t)$, which extends uniquely to 
a $C^\infty$ tangential field $\tau : (\R, 0) \to \R^{1+p}$ of form 
$$
\tau(t) = {\mbox{\rm transpose of}}\left(	1, (a_2/a_1)t^{a_2-a_1} + o(t^{a_2-a_1}), \tau_3(t), 
\dots, \tau_{1+p}(t)\right), 
$$
with $\tau_i(t) = o(t^{a_2-a_1}), (i = 3, \dots, 1+p)$. Therefore $f$ is a frontal. 
Then the tangent map of $f$ is given by $\Tan(f)(t, s) = f(t) + s\tau(t)$. 
By a simple calculation, we have that the Jacobian ideal of $\Tan(f)$ is principal and, in fact, is generated by $st^{a_2-a_1-1}$. Therefore $\Tan(f)$ has a nowhere dense singular locus and turns out to be a frontal 
by Proposition \ref{principal}. 
\qed

\

\noindent
{\it Proof of Theorem \ref{parallel-of-tangent}.} \ 
Let $f : (\R, 0) \to \R^{1+p}$ be a frontal curve and $F = \Tan(f) : 
(\R, 0) \times \R \to \R^{1+p}$ the tangent map of $f$. 
Let $\{ \tau, \mu, \nu_1, \dots, \nu_{p-1}\}$ be an orthonormal frame of $(\Tan(f))^*T\R^{1+p}$ along $f$. 
Here we demand that $\tau$ is a unit section of $T_f$, $\{ \tau, \mu\}$ is an orthonormal frame 
of $T_{\Tan(f)}\vert_{(\R, 0)\times 0}$ and $\{ \nu_1, \dots, \nu_{p-1}\}$ is a parallel orthonormal 
frame of $N_{\Tan(f)}\vert_{(\R, 0)\times 0}$.  Then we extend 
$\{ \tau, \mu, \nu_1, \dots, \nu_{p-1}\}$ trivially to the frame of $(\Tan(f))^*T\R^{1+p}$ 
by the constancy of $T_{\Tan(f)}\vert_{t \times \R}$. 
Then we have $f'(t) = a(t) \tau(t)$ and 
$$
\left(
\begin{array}{c}
\tau'(t) \\
\mu'(t) \\
\nu'_1(t) \\
\vdots \\
\nu'_{p-1}(t)
\end{array}
\right)
= 
\left(
\begin{array}{ccccc}
0 & \kappa(t) & 0 & \cdots & 0 \\
- \kappa(t) & 0 & \ell_1(t) & \cdots & \ell_{p-1}(t) \\
0 & -\ell_1(t) & 0 & \cdots & 0 \\
\vdots & \vdots & \vdots & \ddots & \vdots \\
0 & -\ell_{p-1}(t) & 0 & \cdots & 0
\end{array}
\right)
\left(
\begin{array}{c}
\tau(t) \\
\mu(t) \\
\nu_1(t) \\
\vdots \\
\nu_{p-1}(t)
\end{array}
\right)
$$ 
for uniquely determined functions $a(t), k(t), \ell_1(t), \dots, \ell_{p-1}(t)$. 

Let $\nu$ be a normally parallel, not necessarily unit, normal field along $\Tan(f)$. Then 
$\nu = \sum_{i=1}^{p-1} u_i\nu_i$ for some $u_i \in \R\ (i = 1, \dots, p-1)$. 
Consider parallels 
$$
P(t, s) := {\mbox{\rm{Pal}}}_\nu(\Tan(f))(t, s) = F + \nu = f(t) + s\tau(t) + \nu(t), 
$$ 
of $F = \Tan(f)$ along $\nu$. 
Then the Jacobi matrix of $P$ is given by 
$$
\left(f'(t) + s\tau'(t) + \sum_{i=1}^{p-1} u_i\nu_i'(t), \tau(t)\right) 
= \left(a(t)\tau(t) + s\kappa(t)\mu(t) + \sum_{i=1}^{p-1} u_i(-\ell_i(t)\mu(t)), \tau(t)\right), 
$$
the rank of which is equal to that of the matrix
$\left((s\kappa(t) - \sum_{i=1}^{p-1} u_i\ell_i(t))\mu(t), \ \tau(t)\right)$. 
Therefore the the singular locus $\Sigma(P)$ of the parallel $P$ is given by 
$\{ (t, s) \in (\R, 0)\times \R \mid s \kappa(t) - \sum_{i=1}^{p-1}u_i\ell_i(t) = 0\}$. 
By the assumption, $\kappa(t)$ does not vanish. Therefore  
$\Sigma(P)$ is parametrised by $s = \sum_{i=1}^{p-1}u_i\ell_i(t)/\kappa(t)$. 
Set 
$$
g(t) := f(t) + \sum_{i=1}^{p-1}u_i\frac{\ell_i(t)}{\kappa(t)}\tau(t) + \nu(t).
$$
Then 
$$ 
\begin{array}{rcl}
g'(t) & = & f'(t) + \sum_{i=1}^{p-1}u_i\left(\frac{\ell_i(t)}{\kappa(t)}\tau(t)\right)' + \nu'(t) 
\\
& = & a(t)\tau(t) + \sum_{i=1}^{p-1}u_i\left(\frac{\ell_i(t)}{\kappa(t)}\right)'\tau(t) 
+ \sum_{i=1}^{p-1}u_i\left(\frac{\ell_i(t)}{\kappa(t)}\right)\kappa(t)\mu(t) + \sum_{i=1}^{p-1}u_i(-\ell_i(t)\mu(t))
\\
& = & \left( a(t) + \sum_{i=1}^{p-1} u_i\left\{ \frac{\ell_i(t)}{\kappa(t)}\right\}'\right)\tau(t).
\end{array}
$$
Thus $g$ has the same tangent frame $\tau$ with $f$. The tangent map of $g$ is given by 
$$
\Tan(g)(t, s) = f(t) + \sum_{i=1}^{p-2}u_i\frac{\ell_i(t)}{\kappa(t)}\tau(t) + \nu(t) + s\tau(t) 
= f(t) + \left(s + \sum_{i=1}^{p-2}u_i\frac{\ell_i(t)}{\kappa(t)}\right)\tau(t) + \nu(t). 
$$
Then by the diffeomorphism $(t. s) \mapsto (t, s + \sum_{i=1}^{p-2}u_i\frac{\ell_i(t)}{\kappa(t)})$ 
on $(\R, 0)\times \R$, we have that ${\mbox{\rm{Pal}}}_\nu(\Tan(f))$ is right equivalent 
to $\Tan(g)$. 
\qed

\

\

\noindent 
ISHIKAWA Goo, \\ 
Department of Mathematics, Hokkaido University, \\
Sapporo 060-0810, JAPAN. \\ 
e-mail : ishikawa@math.sci.hokudai.ac.jp 

\end{document}